\documentclass[12pt]{amsart}
\usepackage{txfonts}

\usepackage{amssymb}
\usepackage{mathrsfs}
\usepackage{amsmath}
\usepackage{amssymb,latexsym}

\def\F{\mathbb F} 
\def\Z{\mathbb Z} 

\newtheorem{thm}{Theorem}[section]
\newtheorem{cor}[thm]{Corollary}
\newtheorem{lemma}[thm]{Lemma}

\theoremstyle{definition}
\newtheorem{defn}[thm]{Definition}
\newtheorem{rem}[thm]{Remark}

\numberwithin{equation}{section}

\begin{document}
\baselineskip=17pt

\title [The genus fields of Artin-Schreier extensions ]{The genus fields of Artin-Schreier extensions}

\author[S.Hu and Y.Li]{Su Hu and Yan Li}

\address{Department of Mathematical Sciences, Tsinghua University, Beijing 100084,
China}
\email{hus04@mails.tsinghua.edu.cn
\\liyan\_00@mails.tsinghua.edu.cn}

\begin{abstract} Let $q$ be a power of a prime number $p$. Let $k=\mathbb{F}_{q}(t)$
be the rational function field with constant field $\mathbb{F}_{q}$.
Let $K=k(\alpha)$ be an Artin-Schreier extension of $k$. In this
paper, we explicitly describe the ambiguous ideal classes and the
genus field of $K$ . Using these results we study the $p$-part of the
ideal class group of the integral closure of $\mathbb{F}_{q}[t]$ in
$K$. And we also give an analogy of R$\acute{e}$dei-Reichardt's
formulae for $K$.
\end{abstract}

 \maketitle

\section{Introduction}
In 1951, Hasse~\cite{hasse2} introduced genus theory for quadratic
number fields which is very important for studying the ideal class
groups of quadratic number fields. Later,
Fr$\ddot{o}$hlich~\cite{Frohlich} generalized this theory to
arbitrary number fields. In 1996, S.Bae and J.K.Koo~\cite{Bae}
defined the genus field for global function fields and developed the
analogue of the classical genus theory. In 2000, Guohua
Peng~\cite{Peng} explicitly described the genus theory for Kummer
function fields.

The genus theory for function fields is also very important for
studying the ideal class groups of function fields. Let $l$ be a
prime number and $K$ be a cyclic extension of degree $l$ of the
rational function field $\mathbb{F}_{q}(t)$ over a finite field of
characteristic $\not=l$. In 2004, Wittmann~\cite{Wittmann}
generalized Guohua Peng's results to the case $l\nmid q-1$ and used
it to studied the $l$ part of the ideal class group of the integral
closure of $\mathbb{F}_{q}[t]$ in $K$ following an ideal of
Gras~\cite{Gras}.

Let $q$ be a power of a prime number $p$. Let $k=\mathbb{F}_{q}(t)$
be the rational function field with constant field $\mathbb{F}_{q}$.
Assume that the polynomial $T^{p}-T-D\in k(T)$ is irreducible. Let
$K=k(\alpha)$ with $\alpha^{p}-\alpha=D$. $K$ is called an
Artin-Schreier extension of $k$ (See ~\cite{Hasse}). It is well
known that every cyclic extension of $\mathbb{F}_{q}(t)$ of degree
$p$ is an Artin-Schreier extension. In this paper, we explicitly
describe the genus field of $K$. Using this result we also study the
$p$-part of the ideal class group of the integral closure of
$\mathbb{F}_{q}[t]$ in $K$. Our results combined with
Wittmann~\cite{Wittmann}'s results give the complete results for
genus theory of cyclic extensions of prime degree over rational
function fields.

Let $O_{K}$ be the integral closure of $\mathbb{F}_{q}[t]$ in $K$.
Let $Cl(K)$ be the ideal class group of the Dedekind domain $O_{K}$.
Let $G(K)$ be the genus field of $K$. Our paper is organized as
follows. In Section 2, we recall the arithmetic of Artin-Schreier
extensions. In Section 3, we recall the definition of $G(K)$ and
compute the ambiguous ideal classes of $Cl(K)$ using cohomological
methods. As a corollary, we obtain the the order of Gal($G(K)/K$).
In Section 4, we described explicitly $G(K)$. In Section 5, we study
the $p$-part of $Cl(K)$. And we also give an analogy of
R$\acute{e}$dei-Reichardt's formulae~\cite{RR} for $K$.

\section{The arithmetic of Artin-Schreier extensions}
Let $q$ be a power of a prime number $p$. Let $k=\mathbb{F}_{q}(t)$
be the rational function field. Let $K/k$ be a cyclic extension of
degree $p$. Then  $K/k$ is an Artin-Schreier extension, that is,
$K=k(\alpha)$, where $\alpha^{p}-\alpha=D,\ D\in\mathbb{F}_{q}(t)$
and $D$ can not be written as $x^{p}-x$ for any $x\in k$.
Conversely, for any $D\in\mathbb{F}_{q}(t)$ and $D$ can not be
written as $x^{p}-x$ for any $x\in k$, $k(\alpha)/k$ is a cyclic
extension of degree {p}, where $\alpha^{p}-\alpha=D$. Two
Artin-Schreier extensions $k(\alpha)$ and $k(\beta)$ such that
$\alpha^p-\alpha=D$ and $\beta^p-\beta=D'$ are equal if and only if
they satisfy the following relations,
\begin{equation*}
  \begin{aligned}
          \alpha &\rightarrow x\alpha+B_{0}=\beta,\\
          D&\rightarrow xD+(B_{0}^{p}-B_{0})=D',\\
          x&\in\mathbb{F}_{p}^{*}, B_{0}\in k.
                           \end{aligned}
                           \end{equation*}
(See ~\cite{Hasse} or Artin~\cite{Artin} p.180-181 and p.203-206)
Thus we can normalize $D$ to satisfy the following conditions,

          $$D=\sum_{i=1}^{m}\frac{Q_{i}}{P_{i}^{e_{i}}} + f(t),$$
                  $$(P_{i},Q_{i})=1,\
\ \textrm{and}\ p\nmid e_i,\ \textrm{for}\ 1\leq i\leq m,$$
         $$p\nmid \textrm{deg}(f(t)),\ \textrm{if}\ f(t)\not\in \mathbb{F}_{q},$$
where $P_{i}(1\leq i\leq m)$ are monic irreducible polynomials in
$\mathbb{F}_{q}[t]$ and $Q_{i}(1\leq i\leq m)$ are polynomials in
$\mathbb{F}_{q}[t]$ such that $\textrm{deg} (Q_{i}) < \textrm{deg}
(P_i^{e_i})$. In the rest of this paper, we always assume $D$ has
the above normalized forms and denote
$\frac{Q_{i}}{P_{i}^{e_{i}}}=D_i$, for $1\leq i\leq m$. The infinite
place $(1/t)$ is splitting, inertial, or ramified in $K$
respectively when $f(t)=0$; $f(t)$ is a constant and the equation
$x^{p}-x=f(t)$ has no solutions in $\mathbb{F}_{q}$; $f(t)$ is not a
constant. Then the field $K$ is called real, inertial imaginary, or
ramified imaginary respectively. The finite places of $k$ which are
ramified in $K$ are $P_{1},\cdots,P_{m}$ (p.39 of~\cite{Hasse}).Let $\mathfrak{P}_{i}$  be the place of $K$ lying above $P_{i} (1\leq i\leq m)$.

 Let $P$ be
a finite place of $k$ which is unramified in $K$. Let $(P,K/k)$ be
the Artin symbol at $P$. Then
$$(P,K/k)\alpha=\alpha+\{\frac{D}{P}\}$$  and the Hasse symbol
$\{\frac{D}{P}\}$ is determined by the following equalities:
\begin{equation*}\begin{aligned}
\{\frac{D}{P}\} &\equiv D+D^{p}+\cdots D^{N(P)/p}\textrm{mod} ~P\\ &\equiv (D+D^{p}+\cdots D^{N(P)/p}) \\
&+ (D+D^{q}+\cdots D^{N(P)/q})^{p}\\
&+ \cdots\\
&+  (D+D^{q}+\cdots D^{N(P)/q})^{q/p} \textrm{mod}
~P,\\\{\frac{D}{P}\} &=
\textrm{tr}_{\mathbb{F}_{q}/\mathbb{F}_{p}}\textrm{tr}_{(O_{K}/P)/\mathbb{F}_{q}}(D)
\ \rm{mod}\
  P\end{aligned}
\end{equation*}
(p.40 of ~\cite{Hasse}).

\section{Ambiguous ideal classes }  From this point, we will
use the following notations:
\begin{equation*}
\begin{aligned}
q   ~~~&- ~\textrm{power of a prime number} ~p.\\
k   ~~~&- ~\textrm{the rational function field}~\mathbb{F}_{q}(t). \\
K   ~~~&- ~\textrm{an Artin-Schreier extension of}~k~\textrm{of degree}~p. \\
G   ~~~&- ~\textrm{the Galois group}~Gal(K/k).\\
\sigma  ~~~&- ~\textrm{the geneartor of }~Gal(K/k).\\
S ~~~&- ~\textrm{the set of infinite places of }~K\ \textrm{(i.e,\ the\ primes\ above\ }(1/t)).\\
O_{K}  ~~~&- ~\textrm{the integral closure of}~\mathbb{F}_{q}[t]~\textrm{in}~K. \\
I(K)  ~~~&- ~\textrm{the group of fractional ideals of }~O_{K}.\\
P(K)  ~~~&- ~\textrm{the group of principal fractional ideals of }~O_{K}.\\
P(k)  ~~~&- ~\textrm{the subgroup of}~ P(K)\textrm{ generated by nonzero elements of}~ \mathbb{F}_{q}(t). \\
Cl(K)   ~~~&- ~\textrm{the ideal class group of
 }~O_{K}.\\
H(K)  ~~~&- ~\textrm{the Hibert class field of }~K.\\
G(K)   ~~~&- ~\textrm{the genus field of }~K.\\
U_{K}  ~~~&- ~\textrm{the unit group of }~O_{K}.\\
\end{aligned}
\end{equation*}

\begin{defn}(Rosen~\cite{Rosen})The Hilbert class field $H(K)$ of $K$ (relative to $S$) is the maximal unramified abelian
extension of $K$ such that every infinite places (i.e. $\in\ S$) of
$K$ split completely in $H(K)$.
\end{defn}
\begin{defn}(Bae and Koo ~\cite{Bae}) The genus field $G(K)$ of $K$
is the maximal abelian extension of $K$ in $H(K)$ which is the
composite of $K$ and some abelian extension of $k$.
\end{defn}
For any $G$-module $M$, let $M^{G}$ be the $G$-module of elements of
$M$ fixed by the action of $G$. Without lost of generality, we will
assume $K/k$ is a geometric extension in the rest of this paper. We
have the following Theorem.
\begin{thm} \label{th:1}The ambiguous ideal classes $Cl(K)^{G}$ is a
vector space over $\F_p$ generated by by
$[\mathfrak{P}_{1}],[\mathfrak{P}_{2}],\cdots,[\mathfrak{P}_{m}]$ with dimension
\begin{equation*}dim_{\mathbb{F}_{p}}Cl(K)^{G}=
\begin{cases}
 m-1 &K~\textrm{is real}. \\
  m  &K~\textrm{is imaginary}.
\end{cases}
\end{equation*}

\end{thm}
 Before the proof of the above
theorem, we need some lemmas.
\begin{lemma}\label{nt1}
$H^{1}(G,P(K))=1.$
\end{lemma}
\begin{proof}
From the following exact sequence $$1\longrightarrow
U_{K}\longrightarrow K^{*}\longrightarrow P(K)\longrightarrow 1,$$we
have  $$1\longrightarrow H^{1}(G,P(K))\longrightarrow
H^{2}(G,U_{K})\longrightarrow H^{2}(G,K^{*})\longrightarrow\cdots$$
This is because $K/k$ is a cyclic extension and $H^{1}(G,K^{*})=1$
(Hilbert Theorem 90). Since
\begin{equation}\label{eq:1}
H^{2}(G,U_{K})=\frac{U_{K}^{G}}{NU_{K}}=\frac{\mathbb{F}_{q}^{*}}{(\mathbb{F}_{q}^{*})^{p}}=1,\end{equation}
we have $H^{1}(G,P(K))=1$.
\end{proof}

\begin{lemma}\label{nt2}
If $K$ is imaginary, then $H^{1}(G,U_{K})=1.$
\end{lemma}
\begin{proof}
  Since $U_{K}=\mathbb{F}_{q}^{*}$, we have
$$H^{1}(G,\mathbb{F}_{q}^{*})=\frac{\{x\in\mathbb{F}_{q}^{*}|x^{p}=1\}}{\{x^{\sigma-1}|x\in\mathbb{F}_{q}^{*}\}}=1.$$
\end{proof}

\begin{lemma}\label{nt3}
If $K$ is real, then $H^{1}(G,U_{K})\cong\mathbb{F}_{p}.$
\end{lemma}
\begin{proof} We denote by $\mathscr{D}$
the group of divisors of $K$, by $\mathscr{P}$ the subgroup of
principal divisors. We define $\mathscr{D}(S)$ to be the subgroup of
$\mathscr{D}$ generated by the primes in $S$ and $\mathscr{D}^0(S)$
to be the degree zero divisors of $\mathscr{D}(S)$. From Proposition
14.1 of ~\cite{Rosen2}, we have the following exact sequence
$$(0)\longrightarrow \mathbb{F}_{q}^{*}\longrightarrow U_{K}\longrightarrow
\mathscr{D}^0(S)\longrightarrow Reg\longrightarrow(0),$$ where the
map from $U_{K}$ to $\mathscr{D}^0(S)$ is given by taken an element
of $U_{K}$ to its divisor and $Reg$ is a finite group (See
Proposition 14.1 and Lemma 14.3 of ~\cite{Rosen2}). By Proposition 7
and Proposition 8 of ~\cite{Serre} (p.134), we have
$h(U_{K})=h(\mathscr{D}^0(S))$, where $h(*)$ is the Herbrand
Quotient of $*$. By Equation (\ref{eq:1}), we have
$H^{2}(G,U_{K})=1$. Thus, we can prove this Lemma by showing
$h(\mathscr{D}^0(S))=1/p$.

Let $\infty$ be any infinite place in $S$. Thus $\mathscr{D}^{0}(S)$
is the free abelian group generated by
$(\sigma-1)\infty,(\sigma^{2}-\sigma)\infty,\cdots,(\sigma^{p-1}-\sigma^{p-2})\infty.$
And we have
\begin{equation}\label{eq:2}\mathscr{D}^0(S)=\mathbb{Z}[G](\sigma-1)\infty\cong\frac{\mathbb{Z}[G]}{(1+\sigma+\cdots\sigma^{p-1})}.\end{equation}
Let $\zetaup_{p}$ be a $p$-th root of unity. We have
\begin{equation}\label{eq:3}\frac{\mathbb{Z}[G]}{(1+\sigma+\cdots\sigma^{p-1})}\cong\Z[\zetaup_{p}],\end{equation}
and the above map is given by taken $\sigma$ to $\zetaup_{p}$.  From
(\ref{eq:2}) and (\ref{eq:3}), we have
\begin{equation*}\begin{aligned} H^{1}(G,\mathscr{D}^0(S))&=\frac{ker N
\mathscr{D}^0(S)}{(\sigma-1)\mathscr{D}^0(S)}= \frac{
\mathscr{D}^0(S)}{(\sigma-1)\mathscr{D}^0(S)}\\&\cong
\frac{\frac{\mathbb{Z}[G]}{(1+\sigma+\cdots\sigma^{p-1})}}{(\sigma-1)\frac{\mathbb{Z}[G]}{(1+\sigma+\cdots\sigma^{p-1})}}\cong\frac{\mathbb{Z}[\zetaup_{p}]}
{(\zetaup_{p}-1)}\cong\mathbb{F}_{p}\end{aligned}\end{equation*} and
\begin{equation*}\begin{aligned} H^{2}(G,\mathscr{D}^0(S))&=\frac{
\mathscr{D}^0(S)^G}{N\mathscr{D}^0(S)}=0.\end{aligned}\end{equation*}
Thus $h(\mathscr{D}^0(S))=1/p$.
\end{proof}Proof of Theorem~\ref{th:1}: From the following exact sequence
$$1\longrightarrow P(K)\longrightarrow I(K)\longrightarrow
Cl(K)\longrightarrow 1,$$we have $$1\longrightarrow
P(K)^{G}\longrightarrow I(K)^{G}\longrightarrow
Cl(K)^{G}\longrightarrow H^{1}(G,P(K))\longrightarrow \cdots$$ Since
$H^{1}(G,P(K))=1$ by Lemma~\ref{nt1}, ~we have $$1\longrightarrow
P(K)^{G}\longrightarrow I(K)^{G}\longrightarrow
Cl(K)^{G}\longrightarrow 1.$$ Thus
\begin{equation}\label{eq:4} 1\longrightarrow \frac{P(K)^{G}}{P(k)}\longrightarrow
\frac{I(K)^{G}}{P(k)}\longrightarrow Cl(K)^{G}\longrightarrow
1.\end{equation} From the following exact sequence
$$1\longrightarrow U_{K}\longrightarrow K^{*}\longrightarrow
P(K)\longrightarrow 1,$$ we have
$$1\longrightarrow\mathbb{F}_{q}^{*}\longrightarrow k^{*}\longrightarrow
P(K)^{G}\longrightarrow H^{1}(G,U_{K})\longrightarrow 1$$ and
\begin{equation}\label{eq:5}H^{1}(G,U_{K})\cong\frac{P(K)^{G}}{P(k)}.\end{equation}
Since $\frac{I(K)^{G}}{P(k)}$ is a vector space over
$\mathbb{F}_{p}$ with basis $[\mathfrak{P}_{1}],[\mathfrak{P}_{2}],\cdots,[\mathfrak{P}_{m}]$, by
(\ref{eq:4}), (\ref{eq:5}), Lemma~\ref{nt2} and Lemma~\ref{nt3}, we
get the desired result.
\begin{rem}If $K$ is real, it is an interesting question to find explicitly the
relation satisfied by $[\mathfrak{P}_{1}],[\mathfrak{P}_{2}],\cdots,[\mathfrak{P}_{m}]$ in $Cl(K)^G$.
By Lemma \ref{eq:5}, if we can find a nontrivial element $\bar{u}$
of $H^{1}(G,U_{K})$, then by Hibert 90, we have $u=x^{\sigma-1}$,
where $u\in U_{K}$ and $x\in K$. It is easy to see that
$$\sum_{i=1}^{m}ord_{\mathfrak{P}_i}(x)[\mathfrak{P}_i]=0$$
in $Cl(K)^G$.
\end{rem}
From Proposition 2.4 of ~\cite{Bae}, we have
\begin{equation}\label{eq:6}
Gal(G(K)/K)\cong Cl(K)/(\sigma-1)Cl(K)\cong Cl(K)^{G}.
\end{equation}
(It should be noted that the last isomorphism is merely an
isomorphism of abelian groups but not canonical). Therefore, we get
\begin{cor}\label{cor}
\begin{equation*}{\#Gal(G(K)/K)=}
\begin{cases}
p^{m-1} &K~\textrm{is real}. \\
p^{m} &K~\textrm{is imaginary}.
\end{cases}
\end{equation*}
\end{cor}
\section{The genus field $G(K)$}
In this section, we prove the following theorem which is the main
result of this paper.

\begin{thm}\label{main}\begin{equation*}{G(K)=}\begin{cases}
 k(\alpha_{1},\alpha_{2},\cdots,\alpha_{m})&K~\textrm{is real}. \\
      k(\beta,\alpha_{1},\alpha_{2},\cdots,\alpha_{m})&K~\textrm{is imaginary}.
\end{cases}
\end{equation*}
Where
$\alpha_{i}^{p}-\alpha_{i}=D_i=\frac{Q_{i}}{P_{i}^{e_{i}}}(1\leq
i\leq m),\beta^{p}-\beta=f(t),$ and $D_i,Q_{i},P_{i},f(t)$ are
defined in Section 2.
\end{thm}We only prove the imaginary case. The proof is the same for
the real case. Since
$$(\sum_{i=1}^{m}\alpha_{i}+\beta)^{p}-(\sum_{i=1}^{m}\alpha_{i}+\beta)=\sum_{i=1}^{m}\frac{Q_i}{P_{i}^{e_{i}}}+f(t)=D,$$
we can assume $\alpha=\sum_{i=1}^{m}\alpha_{i}+\beta.$ ~Before the
proof of the above theorem, we need two lemmas.
\begin{lemma}\label{nt4}$E=k(\beta,\alpha_{1},\alpha_{2},\cdots,\alpha_{m})$ is
an unramified abelian extension of $K$.
\end{lemma}
\begin{proof}
Let $P$ be a place of $k$ and let $(1/t)$ be the infinite place of
$k$. If $P\not=P_{1},P_{2},\cdots,P_{m},(1/t)$, then $P$ is
unramified in $ k(\beta),k(\alpha_{i})(1\leq i\leq m)$, hence
unramified in $E$. Otherwise, without lost of generality, we can
suppose $P=P_{1}$. Since $\alpha=\sum_{i=1}^{m}\alpha_{i}+\beta,$ we
have $E=Kk( \alpha_{2},\cdots,\alpha_{m},\beta)$. Thus $P=P_{1}$ is
unramified in $k(\alpha_{2},\cdots,\alpha_{m},\beta)$, hence
unramified in $E/K$.
\end{proof}
\begin{lemma}\label{nt5}The infinite places of $K$ are split completely in $E=k(\beta,\alpha_{1},\alpha_{2},\cdots,\alpha_{m})$.
\end{lemma}
\begin{proof}
Since $\alpha=\sum_{i=1}^{m}\alpha_{i}+\beta,$ we have
$E=Kk(\alpha_{1}, \alpha_{2},\cdots,\alpha_{m})$. Since the infinite
place $(1/t)$ of $k$ splits completely in $k(\alpha_{1},
\alpha_{2},\cdots,\alpha_{m})$, hence $(1/t)$ also splits completely
in $E/K$.
\end{proof}
Proof of Theorem~\ref{main}:From Lemma~\ref{nt4} and~\ref{nt5},~we
have
\begin{equation}\label{eq:8}k(\alpha_{1},\alpha_{2},\cdots,\alpha_{m},\beta)\subset
G(K).\end{equation}  Comparing ramifications, $
k(\beta),k(\alpha_{i})(1\leq i\leq m)$ are linearly disjoint over
$k$, so
$$[k(\alpha_{1},\alpha_{2},\cdots,\alpha_{m},\beta):k]=p^{m+1}$$ and
$$[k(\alpha_{1},\alpha_{2},\cdots,\alpha_{m},\beta):K]=p^{m}.$$ Thus from
Corollary~\ref{cor} and (\ref{eq:8}), we get the result.

\section{The $p$-part of $Cl(K)$}
If $l$ is a prime number, $K$ is a cyclic extension of $k$ of degree
$l$, and $\mathbb{Z}_{l}$ is the ring of $l$-adic integers, then
$Cl(K)_{l}$ is a finite module over the discrete valuation ring
$\mathbb{Z}_{l}[\sigma]/(1+\sigma+\cdots+\sigma^{l-1})$. Thus its
Galois module structure is given by the dimensions:
$$\lambda_{i}=dim(Cl(K)_{l}^{(\sigma-1)^{i-1}}/Cl(K)_{l}^{(\sigma-1)^{i}})$$
for $i \geq 1$, these quotients being $\mathbb{F}_{l}$ vector spaces
in a natural way. In number field situations, the dimensions
$\lambda_{i}$ have been investigated by R$\acute{e}$dei~\cite{RR}
for $l=2$ and Gras~\cite{Gras} for arbitrary $l$. In function field
situations, these dimensions $\lambda_{i}$ have been investigated by
Wittmann for $l\not=p$. In this section, we give a formulae to
compute $\lambda_{2}$ for $l=p$. This is an analogy of
R$\acute{e}$dei-Reichardt's formulae~\cite{RR} for Artin-Schreier
extensions.

If $K$ is imaginary, as in the proof of Theorem~\ref{main}, we
suppose that $K=k(\alpha)$, where
$\alpha=\sum_{i=1}^{m}\alpha_{i}+\beta$.  We have the following
sequence of maps \begin{equation*}
\begin{aligned} Cl(K)^{G}
\longrightarrow Cl(K)/(\sigma-1)Cl(K)\cong
Gal(G(K)/K)&\hookrightarrow Gal(G(K)/k)\\\cong Gal
(k(\alpha_{1})/k)&\times\cdots\times Gal (k(\alpha_{m})/k)\times
Gal(k(\beta)/k).\end{aligned}\end{equation*}  Considering $[\mathfrak{P}_{i}]\in
Cl(K)^{G}~(1\leq i\leq m)$  under these maps, we have
\begin{equation*}
\begin{aligned}
[ ] [\mathfrak{P_{i}}] \longmapsto \bar{[\mathfrak{P}_{i}]} &\longmapsto (\mathfrak{P}_{i},G(K)/K) \longmapsto
(\mathfrak{P}_{i},G(K)/K)\\&\longmapsto((P_{i},
 k(\alpha_{1})/k),\cdots,(P_{i},
k(\alpha_{m})/k),(P_{i},k(\beta)/k)),\end{aligned}\end{equation*}
where the $i$-th component is $(\mathfrak{P}_i, G(K)/K)|_{k(\alpha_i)}$.

We define the R$\acute{e}$dei matrix $R=(R_{ij})\in M_{m\times
m}(\mathbb{F}_{p})$ as following:
$$R_{ij}=\{\frac{D_{j}}{P_{i}}\},~\textrm{for}~ 1\leq i,j\leq m,i\not=j,$$
and $R_{ii}$ is defined to satisfy the equality:
$$\sum_{j=1}^{m}R_{ij}+\{\frac{f}{P_{i}}\}=0.$$
From the discussions in section 2, we have
\begin{equation*}\begin{aligned}(\mathfrak{P}_{i},G(K)/K)\alpha &=\alpha
,\\(\mathfrak{P}_{i},G(K)/K)\alpha_{j} &=\alpha_{j}+\{\frac{D_{j}}{P_{i}}\},
for\ i\neq j\\(\mathfrak{P}_{i},G(K)/K)\beta
&=\beta+\{\frac{f}{P_{i}}\},\end{aligned}\end{equation*} so it is
easy to see the image of $Cl(K)^{G} \rightarrow
Cl(K)/(\sigma-1)Cl(K)$ is isomorphic to the vector space generated
by the row vectors
$(R_{i1},R_{i2},\cdots,R_{im},\{\frac{f}{P}_{i}\})\ (1\leq i\leq
m).$

We conclude that
\begin{equation*}\begin{aligned}\lambda_{2}&=dim_{\mathbb{F}_{p}}(Cl(K)_{l}^{(\sigma-1)}/Cl(K)_{l}^{(\sigma-1)^{2}})=dim_{\mathbb{F}_{p}}(Cl(K)^{(\sigma-1)}/Cl(K)^{(\sigma-1)^{2}})\\&=dim_{\mathbb{F}_{p}} \ker ( Cl(K)^{G}
\rightarrow
Cl(K)/(\sigma-1)Cl(K))\\&=dim_{\mathbb{F}_{p}}Cl(K)^{G}-dim_{\mathbb{F}_{p}}
Im ( Cl(K)^{G} \rightarrow Cl(K)/(\sigma-1)Cl(K))\\&= m-rank(R).
\end{aligned}\end{equation*}

Since the proof of real case is similar, we only give the results
and sketch the proof.

If $K$ is real, from the discussions in section 2, we have $f(t)=0$,
so $$D=\sum_{i=1}^{m}D_{i}.$$ We define the R$\acute{e}$dei matrix
$R=(R_{ij})\in M_{m\times m}(\mathbb{F}_{p})$ as following:
$$R_{ij}=\{\frac{D_{j}}{P_{i}}\},~\textrm{for}~ 1\leq i,j\leq m,i\not=j,$$
and $R_{ii}$ is defined to satisfy the equality:
$$\sum_{j=1}^{m}R_{ij}=0.$$
The same procedure as the imaginary case shows that the image of
$Cl(K)^{G} \rightarrow Cl(K)/(\sigma-1)Cl(K)$ is isomorphic to the
vector spaces generated by the row vectors of R$\acute{e}$dei
matrix. Thus

\begin{equation*}\begin{aligned}\lambda_{2}&=dim_{\mathbb{F}_{p}}Cl(K)^{G}-dim_{\mathbb{F}_{p}}
Im ( Cl(K)^{G} \rightarrow Cl(K)/(\sigma-1)Cl(K))\\&= m-1-rank(R).
\end{aligned}\end{equation*}
\begin{thm}If $K$ is imaginary, then $\lambda_2=m-rank(R)$; if $K$
is real, then $\lambda_2=m-1-rank(R)$, where $R$ is the
R$\acute{e}$dei matrix defined above.

\end{thm}
If $p=2$, then $\sigma$ acting on $Cl(K)$ equal to $-1$. So
$\lambda_1$, $\lambda_2$ equal to the 2-rank, 4-rank of ideal class
group $Cl(K)$, respectively. In particular, the above theorem tells
us the 4-rank of ideal class group $Cl(K)$ which is an analogue of
classical R$\acute{e}$dei-Reichardt's 4-rank formulae for narrow ideal
class group of quadratic number fields.

\end{document}